%% file: Dehn-twist.tex
\documentclass[production,11pt]{amsart}
\usepackage{hyperref,amsmath,amsthm,amsfonts,amscd,flafter,epsf}
\def\sqr#1#2{{\vcenter{\vbox{\hrule height.#2pt
	\hbox{\vrule width.#2pt height#1pt \kern#1pt
	\vrule width.#2pt}
	\hrule height.#2pt}}}}
\def\square{\mathchoice\sqr64\sqr64\sqr{2.1}3\sqr{1.5}3}
\begin{document}
\input{command.tex}
\title{Floer Cohomology of Certain Pseudo-Anosov Maps}
\author{Eaman Eftekhary} \thanks{ Mathematics Department, Princeton University,
Princeton, NJ 08544\\
\ \ \ E-mail: eeftekha@math.princeton.edu}
\begin{abstract}
Floer cohomology is computed for the elements of the mapping class
group of a surface
$\Sigma$ of genus $g>1$ which are compositions of positive and negative Dehn
twists along  loops in $\Sigma$ forming a tree-pattern. The computations cover a certain class of
pseudo-Anosov maps.
\end{abstract}
\maketitle
{\em {Dedicated to S. Shahshahani on the occasion of his 60-th birthday}}
\section{introduction}
The symplectic Floer cohomology (homology) is associated with a compact
symplectic $2n$-manifold
$(M^{2n},\omega)$, where $L_1$ and $L_2$ are two
transversal Lagrangian submanifolds   of dimension $n$ and $\omega$ is the symplectic form.
The Floer complex
is freely generated by the
intersection points of $L_1$ and $L_2$. The boundary maps count the
pseudo-holomorphic disks between these points with the boundary on the
Lagrangian submanifolds (see \cite{Fl}). \\

For any symplectic form $\omega$ on a surface $\Sigma$ of genus $g>1$, there
are diffeomorphisms $f:\Sigma \rightarrow \Sigma$ in each isotopy class in
the mapping class group of $\Sigma$, which preserve $\omega$. Putting the symplectic form
$\omega \times (-\omega)$ on $\Sigma \times \Sigma$, the diagonal $\Delta$ and
the graph $\Gamma_f$ of $f$ will become Lagrangian. One may ask for the
symplectic Floer cohomology of this setting, denoted by $HF^*(f)$. There are
some technical issues that may arise, in order to make this cohomology group
 an invariant of the isotopy class. It is shown in
\cite{Sei2} that the Floer cohomology of an isotopy class in the mapping
class group is a well-defined $\mathbb{Z}/2\mathbb{Z}$ -graded $H^*(\Sigma,\mathbb{Z}/2
\mathbb{Z})$
-module.
Monotonicity is necessary for using a particular representative of a class in
the mapping class group for Floer cohomology computations (see \cite{Sei2}).
In \cite{Sei1}, this Floer cohomology is computed for a combination of
Dehn twists along a disjoint union of curves. \\

Gautschi extended this result to the case of algebraically finite maps (\cite{Gau}).
These are obtained by choosing disjoint loops $\alpha_i$ on the surface and composing
the positive Dehn twists along $\alpha_i$s with finite order automorphisms of
$\Sigma \setminus \cup \alpha_i$ which satisfy  certain boundary conditions.
By Thurston's classification of surface diffeomorphisms these are the maps with no
pseudo-Anosov components.\\

On the other hand the mapping class group is generated by positive Dehn twists
along certain simple closed curves on the surface. This means that each class
has a representative given by composition of positive Dehn twists along a
collection of curves $\alpha_i, i=1,2,...,n$ on $\Sigma$. However, they may have
many intersections. Hence it is important to be able to push the computation of
Seidel to the case of loops with intersections.\\

In this paper, the result of \cite{Sei1} is extended in a
different direction, by allowing the loops to intersect each other in a
``nice'' way. \\

In particular many pseudo-Anosov maps may be obtained by allowing these types of
intersections (for instance,
see example~\ref{ex:dehn.1}). These are the first computations for
pseudo-Anosov diffeomorphisms.
There are no holomorphic orbits between the generators of the complex if the
composition of twists represents a pseudo-Anosov class.
It is interesting to study whether in general this is the case . \\
\\
\begin{defn}
A set of data $(\Sigma ,\alpha_1,\alpha_2,...,\alpha_n)$  will be called
an \emph{acceptable setting} if $\Sigma$ is a surface of genus $g>1$ and
$\alpha_1,...,\alpha_n$ are simple closed loops  on $\Sigma$
 satisfying the following conditions: Any two of the loops intersect
 at most once and transversely. If we form the
\emph{intersection graph} $G$ with vertices $1,2,...,n$, and connect
$i,j$ iff $\alpha_i,\alpha_j$ intersect each other, then $G$ is a
\emph{forest} (i.e. it does not contain any loops).
Furthermore, we assume that no $\alpha_i$ is homologically trivial.
\end{defn}
The map will be the combination of positive Dehn twists along $\alpha_i$s.\\

The content of section~\ref{sec:2} will be the construction of some
appropriate Morse function on $\Sigma$. Its Hamiltonian flow  will be composed with
the combination of twists, in order to make the diagonal
transversal to the graph of the function in $\Sigma \times \Sigma$. The fixed points
of this composition will be in $1-1$ correspondence
with the generators of the Floer complex. \\
In section~\ref{sec:3} we will study the moduli spaces of $J$-holomorphic
disks between the fixed points. In particular we will prove some energy bounds.
These bounds will be used in section~\ref{sec:4} to prove the main result
of this paper (theorem~\ref{thm:main}):
\\
\\
 \begin{thm}
Let $(\Sigma,\alpha_1,...,\alpha_n)$ be an acceptable
setting, $C=\cup_i \alpha_i$ and $T$ be the composition of positive Dehn twists
along $\alpha_i$s in some order. Then $HF^*(T)\cong H^*(\Sigma
,C)$ as $H^*(\Sigma,\mathbb{Z}/2\mathbb{Z})$-modules where
$H^*(\Sigma)\cong HF^*(Id)$ acts on the right
hand side by the cup product and on the left hand side by quantum cup
product. (see theorem~\ref{thm:main.general} for a more general
result).
\end{thm}
It is well-known that
in the most familiar case, $HF^*(Id) \cong H^*(\Sigma)$ and the quantum cup
product reduces to the ordinary cup product (see \cite{LT}, \cite{PS} or
\cite{RT}).\\
\\
\\
\textbf{Acknowledgment.} No words can express my thankfulness to Z. Szab\'o,
for his patience in correcting my mistakes, and sharing his insight with me.
Some ideas presented here, are suggested by P. Seidel, who kindly taught me a
lot on this problem. Many thanks go to
M. Hutchings, P. Ozsv\'ath and M. Ajoodanian for the great conversations that
we had.\\

\section{Appropriate Morse Functions}
\label{sec:2}
 We will construct a Morse function
inductively, and use its Hamiltonian  flow to get a  function with the required
trasversality. We will implicitly assume that the loops $\alpha_i, i=1,...,n$ are
smooth.\\
\\
\begin{prop}
\label{prop:morse}
Suppose that $\alpha_i, i=1,2,...,n$ are simple closed
loops on a surface $\Sigma$. Then it is possible to find a Morse function $h:\Sigma
\rightarrow \mathbb{R}$ such that the curves $\alpha_i$ are unions of the flow line
of $h$ and the followings are satisfied:\\
\\
\textbf{(H1)} $h(p)>4$ for all points $p \in \alpha_i, i=1,...,n$.\\
\textbf{(H2)} If $p$ is a critical point of $h$ and $h(p)>2$ then $p\in \alpha_i$
for some $i$.\\
\textbf{(H3)} If $\{q\}=\alpha_i \cap \alpha_j$, then $q$ is an index-$1$ critical
point for $h$\\
\textbf{(H4)} For a critical point $q \in \alpha_i$ of $h$ which is not an intersection
point with other $\alpha_j$s, either $q$ is an index-$2$ critical point, or
it has index $1$ and $h|_{\alpha_i}$ has a local minimum at $q$.
\end{prop}
\begin{proof} The proof is by induction on the number of loops. For a union of
disjoint loops it is easy to find such a function. Now suppose that a map
$h_0:\Sigma \rightarrow \mathbb{R}$ is constructed for $\alpha_1,...,\alpha_{n-1}$
with the above properties. $\alpha_n$ meets the other curves in some points, say
$q_1,...,q_m$ which we may assume to be different from the critical points of
$h_0$. Suppose that $\{q_i\}=\alpha_n \cap \alpha_{j_i}$ for some $j_i \in \{
1,2,...,n-1\}$. Change $h_0$ near each $q_i$ by introducing a pair of cancelling
critical points on $\alpha_{j_i}$ of indices $2,1$ (see figure~\ref{fig:1}(a))
. Choose them such that $q_i$
is the index-$1$ critical point and $h|_{\alpha_{j_i}}$ has a local minimum at
$q_i$.\\

We may assume that near $q_i$, $\alpha_n$ is the union of the two flow lines going out of this
critical point. Let us denote by $\Phi^T$ the gradient flow of the new Morse function up to
time $T$. If $T$ is chosen big enough, $\beta=\Phi^T(\alpha_n)$ will be isotopic
to $\alpha_n$, and  on $h_0^{-1}[2,\infty)$, $\beta$ is a union of flow lines.
Furthermore, $\beta$ has no intersection with other $\alpha$-curves other than $q_i$.
 $\Sigma_0=h_0^{-1}(-\infty,3]$ is a manifold, which is of the form $(2,3]\times
S^1$ near each boundary component. $h_0$ is the projection onto the second factor on
this neighborhood and $\gamma=\beta \cap \Sigma_0$ will be a union of disjoint paths on
$\Sigma_0$. In this situation the following lemma is fairly easy to prove: \\
\begin{figure}
\mbox{\vbox{\epsfbox{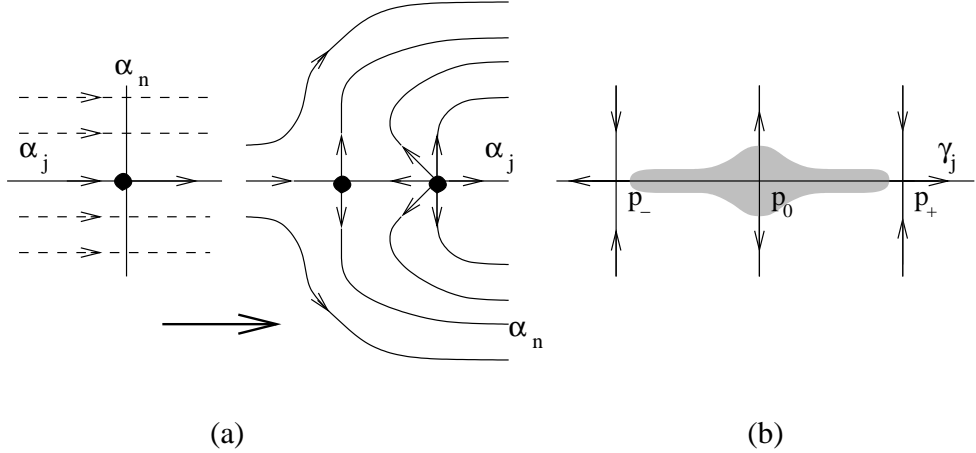}}}
\caption{\label{fig:1}
{(a) Near an intersection point $\alpha_j \cap \alpha_n$ we may introduce
an indext-$2$ critical point on $\alpha_j$ and an index-$1$ critical point
in the intersection, to make $\alpha_n$ locally a flow line. (b) By changing
$h$ in the shaded area, we may assume that $h(p_0)$ has any value less that
$min\{h(p_-),h(p_+)\}$}}
\end{figure}
\\
\begin{lem}
Suppose that $\Sigma_0$ is a surface and $C_1,...,C_l$ are the boundary circles with
tubular neighborhoods $N_i=[2,3]\times C_i \subset \Sigma_0$. Suppose that a union $\gamma$ of
disjoint paths $\gamma_1,...,\gamma_k$ with ends on the boundary circles is of the form
$[2,3]\times \{p_1^i,...,p_{k_i}^i\}$ in each $N_i$. Then there exists a Morse function
$h'$ on $\Sigma_0$ which is the projection onto the first factor on each $N_i$, has value less than
$2$ elsewhere and each $\gamma_i$ is a union of the two flow lines going into an index-$1$
critical points of $h'$.\end{lem}
$\square$\\

Using this lemma we may change $h_0$ on $\Sigma_0$ to be equal to $h'$. Then
$\alpha_n$ is cut into several pieces which are flow lines of the new Morse function
$h'$, say $\alpha_n=
\gamma_1 \cup ...\cup \gamma_k$. Each $\gamma_j$ starts from an index-$1$ critical
point $p_{-}$ and goes down to another index-$1$ point $p_0$, then goes back up to
a third index-$1$ critical point $p_{+}$ (figure~\ref{fig:1}(b)). For these points we have
$h'(p_+),h'(p_-)>4$ and $h(p_0)<2$. By changing $h'$ in the shaded area of
figure~\ref{fig:1}(b) we may move $h'(p_0)$ within the range
$min\{h(p_-),h(p_+)\}>h(p_0)>4$ such that  the flow lines between the critical points
remain the same. This completes the induction.
\end{proof}

 From now on, we will assume that $(\Sigma, \alpha_1,...,\alpha_n)$ is an
acceptable setting and that the positive Dehn twist $T_i$ along $\alpha_i$
is done in small tubular neighborhood of the loop. If we choose the symplectic form
to be standard in the neighborhood of $\alpha_i$s, this map may be assumed to be a
symplectomorphism.\\

\begin{figure}
\mbox{\vbox{\epsfbox{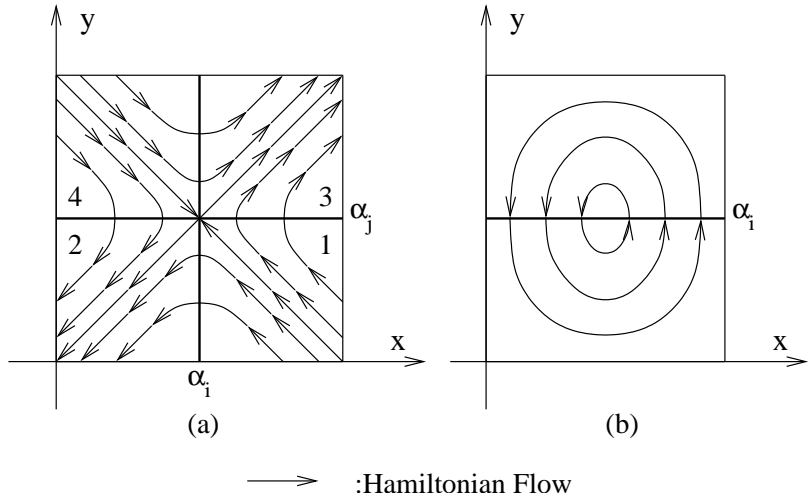}}}
\caption{\label{fig:2}
{Local behavior;
(a) The intersection at index$-1$ critical points of $h$ of $\alpha_i,\alpha_j$
(b) $\alpha_i$ passes through an index$-2$ critical point of $h$}}
\end{figure}

Now it is time to show how one should use this Morse function $h$ to perturb the
combination $T$ of the twists.\\
Let $X_h$ denote the vector field satisfying:
\begin{equation}
\omega (X_h(x),\zeta)=-dh_x\zeta \ \ \ \ \forall \ \ \zeta \in T_x\Sigma,
\end{equation}
where $h$ is the Morse function constructed above. One
may consider the flow of this vector field, denoted by $H^t(x)$, which is
called the \emph{Hamiltonian flow} of $h$.
Put $T^\epsilon :=H^\epsilon \circ T$.
The following result determines the fixed points of $T^
{\epsilon}$:\\
\\
\begin{thm}
\label{thm:fixed}
Suppose that $h$ is the Morse function given by proposition~\ref{prop:morse}.
Then for $\epsilon >0$ small enough, the only fixed
points of
the map $T^\epsilon$ constructed above, are the critical points $p$
of $h$ with $h(p)<2$.\end{thm}
\begin{proof}
First of all notice that if $p$ is a critical point of $h$, with
$h(p)<2$, then it will be a fixed point of $T^\epsilon$. We claim that these
are actually the only fixed points. Away from small strips around the
loops, $T$ is the identity map and the above claim is trivial. Now, assume that
$p$ is a fixed point in this region.
$H^\epsilon$ will not change the ``height'' of the points (i.e. $h(H^\epsilon (x))=h(x)$)
and each point is only moved slightly by $H^\epsilon$. Since $p$ is a fixed point
of $T^\epsilon$, $T$ should have the same property.
This is not the case away from the critical points of $h$.\\

In fact if $p$ is very close to a loop $\alpha_i$, $T_i$ will map $p$ far away from
its initial position. Since the intersection graph of $\alpha_j$s does not have any
loops, it is not possible to compose this twist with other ones and take $p$ back near
its initial position. It is implied then that $p$ should be twisted only slightly
by each of the $T_i$s. Now if $p$ is away from the critical points of $h$, it
will be twisted only by one of the maps $T_i$. This twist will either increase or
decrease the value of $h$ since $\alpha_i$ is a gradient flow line.
So $p$ can not be a fixed point unless it is around one of the critical points of $h$.\\

We should then analyze what happens near critical points of $h$.
Because of the assumption on $h$ the following cases are the only possibilities:\\
\\

\textbf{1)} $p$ is around an index-$1$ critical point $q$ of $h$ where
$\{q\}=\alpha_i \cap \alpha_j$. The two loops are
locally the stable and unstable manifolds of $q$. The local picture is
shown in  figure~\ref{fig:2}(a). One gets $4$ regions
labeled $1,2,3,4$. For the fixed point $p$ of $T^
\epsilon$, the twists and the flow do not change the region of $p$ if $\epsilon$ is
chosen to be small enough.
Let $y(p)$ denote the $y$-coordinate of the point $p$ in the $xy$-plane
shown in figure~\ref{fig:2}(a).  By
analyzing the local behavior of $H^\epsilon$ and $T$, one may find out that
in the regions 1,3, for any point $p$, $y(T(p))>y(p)$ and $y(H^\epsilon(p))>y(p)$.
Similarly in the regions 2,4, $y(H^\epsilon(p)),y(T(p))<y(p)$. So
$T^\epsilon$ may not have a fixed point near such an intersection point $q$.
\\

\textbf{2)} The next possible case is a local minimum for $h|_{\alpha_i}$
at some index one critical point. The only difference from figure~\ref{fig:2}(a)
is that there is no $\alpha_j$ in the picture.
The argument is identical to the above,
since $T,H^\epsilon$ still increase the $y$-coordinate in the regions 1,3 and
decrease it in 2,4.\\
\\
\textbf{3)} An index-$2$ critical point of $h$ on a curve $\alpha_i$;
The local picture is shown in figure~\ref{fig:2}(a). Here there are two regions
and the $x$-coordinate is increased by $T,H^\epsilon$ in one of them and is
decreased in the other one. Again there can not be any fixed points.
\end{proof}

If $h$ is our Morse function, $h^{-1}(3)$ will be a submanifold, hence a
collection of circles on $\Sigma$. One may assume that near these circles
the manifold is a product $[2 ,4] \times S^1$, $h$ is the
projection on the first factor, and that $\omega$ and the almost complex
structure are standard on $[2,4]\times S^1$. Following
an argument similar to that of Gautschi's (\cite{Gau}) we choose a subset
$I_0=[3-\delta,3+\delta]\subset [2,4]=I$ and  consider a bump
function $\lambda_R$ which takes the value $R$ on $I_0$ and is equal to $1$ near the
ends of $I$. Then we may think of $\lambda_R$ as a function on $\Sigma$ which is
equal to the constant function $1$ outside the above product neighborhoods and on
these neighborhoods is given by the value of the above function on the first factor.\\

Multiply the symplectic form by $\lambda_R$, but keep the almost complex structure
fixed. Call the new symplectic form $\omega_R$. It is clear that
the moduli space of holomorphic curves does not change.
 We will use $\omega_R$ for large $R$
to reduce the boundary computations to the case of Hamiltonian isotopy.\\
\\

\section{Energy Bounds}
\label{sec:3}
Before going through the next few technical steps, let us see how they are
related to the bounds on the energy of pseudo-holomorphic disks connecting
the fixed points.
Fix a symplectic form $\omega$ on $\Sigma$ such that
$T$ becomes a symplectomorphism and $\omega$ is standard in the parts described
in the last paragraph of section~\ref{sec:2}.\\

Keeping the notation of section~\ref{sec:2},
for $x_-,x_+\in $Fix$(T^{\epsilon})$ let $\pi_{2}(x_-,x_+)$
be the space of isotopy classes of maps
$\varphi:\mathbb{R} \times \mathbb{R} \rightarrow \Sigma^R$
s.t.
\begin{equation}
\varphi(s,t)=T^{\epsilon}(\varphi(s,t+1)) \ \ \ \  \forall s,t
\end{equation}
\begin{displaymath}
\lim_{s \rightarrow \pm \infty} \varphi(s,t) =  x_{\pm}. \ \ \ \
\end{displaymath}

For $\varphi \in \pi_{2}(x_-,x_+)$ define $\mathcal{M} (\varphi)$
to be
the space of maps $u:\mathbb{R}\times  \mathbb{R} \rightarrow \Sigma$
representing $[\varphi]$ s.t.
\begin{equation}
\frac{du}{ds}+J(u(s,t))\frac{du}{dt}=0 \ \ \ \ \forall s,t.
\end{equation}

Denote the expected dimension of this space
by $\mu (\varphi)$.
$\mathcal{M}(\varphi)$ is obviously equipped with an $\mathbb{R}$-action
which is the translation in the second factor.\\
We consider the Floer complex generated by $
\{p\}_{p\in \textrm{Fix}(T^{\epsilon})}$. The boundary maps are defined by:
\begin{equation}
\partial (p):=
\sum_{\substack{q\in \textrm{Fix} (T^{\epsilon})\\
\varphi \in \pi_{2}(p,q)\\
\mu (\varphi)=1}} \# (\mathcal{M}(\varphi)/\mathbb{R})\cdot q \ .
\end{equation}

If the almost complex structure $J$ and the perturbation are generic enough
then this will define a homology group denoted by $HF^*(T)$.
\\
\\
\begin{prop}
\label{prop:disks}
Suppose that
$\Sigma \setminus \cup_{i=1}^{n} \alpha_i=\cup_{j=1}^{m}V_j$, where $V_j$s are
 the connected components. If $x_-,x_+\in \textrm{Fix}(T^{\epsilon})$ are not
in the same $V_j$, then $\pi_{2}(x_-,x_+)$ will be empty.\end{prop}

As a result we may restrict ourselves to the fixed points on the same piece.
The next result is:\\
\\
\begin{prop}\label{prop:energy} For $x_-,x_+\in \textrm{Fix}(T^{\epsilon})$ in
the same piece $V_j$ and $\varphi \in \pi_{2}(x_-,x_+)$ , $u\in
\mathcal{M}(\varphi)$, the energy of $u$ is $E(u)= n(\varphi)
(\int_{\Sigma} \omega _R) +
\epsilon(h(x_+)-h(x_-))$, for some integer $n(\varphi)\geq 0$ that only
depends on the homotopy class $\varphi$.\end{prop}

Note that the energy depends on the metric and if the metric is defined
using the almost complex structure $J$ and the symplectic form $\omega_R$,
then $\omega_R$ will naturally enter to any formula for the energy.\\
\\

This energy bound will be used later to rule out some potential holomorphic orbits.\\
In the rest of this section we will focus on the proof of the above two
statements.\\
Put
$T_{\pi}:=T_{\pi(n)} \circ T_{\pi(n-1)} \circ ... \circ T_{\pi(1)}$, for any
element $\pi \in S_n$ (the symmetric group in $n$ letters). Here $T_i$ is the
positive Dehn twist along $\alpha_i$. The map in question
will be $T=T_{Id}$. We claim the following:
\\
\\
\begin{lem}\label{lem:1}
Associated with any $\pi \in S_{n}$, is a diffeomorphism
$f:\Sigma \rightarrow \Sigma$ such that $T_{\pi} = f \circ T_{Id} \circ f^{-1}$.
Moreover, $f$ may be chosen to be a combination of the twists $T_{i}$ and
their inverses.\end{lem}
\begin{proof} To prove the lemma, we use the following simple
combinatorial fact:\\
\\
\begin{lem}\label{lem:2}
Let $G$ be a \emph{forest}, with vertices labeled $1,2,...,n$.
Put these numbers around a circle (with some arbitrary order)
. At each step one may interchange
$i,j$ on the circle, if they are not adjacent in $G$, but are neighbors on the
circle. Using these moves, it is then possible to get any permutation
$\pi (1),\pi (2),...,\pi (n)$
of $1,2,...,n$ around the circle in finitely many steps.\end{lem}$\square$\\

If $G$ is the associated graph of  the loops $\alpha_{i}, i=1,...,n$
and $i,j$ are not adjacent in $G$, then
$T_{i} \circ T_{j} = T_{j} \circ T_{i}$.
So for $\pi \in S_n$ with
$\pi(l)=i,\pi(l+1)=j$, $T_{\pi}=T_{\sigma}$ where
$(i,j) \circ \pi = \sigma$. Also, if $\pi(1)=i,\pi(n)=j$ with $i,j$ not
adjacent in $G$ then for $\sigma = (1,n) \circ \pi$ we find
$T_{\sigma} = T_{\pi(1)} \circ T_{\pi(n)}^{-1} \circ T_\pi \circ T_{\pi(1)}^
{-1} \circ T_{\pi(n)} = f_{\pi(1),\pi(n)} \circ T_{\pi} \circ f_{\pi(1),
\pi(n)}^{-1}$
where
$f_{i,j}= T_{i} \circ T_{j}^{-1}$.
In order to get to $T_\pi$ from $T_{Id}$, one may then follow the steps
determined by lemma~\ref{lem:2}.\end{proof}

\begin{lem} \label{lem:3}
Let $T,\alpha_i$ be as above. If $\beta$ is a path on $\Sigma$ such that
$T$ fixes the ends of $\beta$ and
$T\beta$ is homologous to $\beta$, then $\beta$ has zero
intersection number with all $[\alpha_i]$s.\end{lem}

Here by ``$T\beta$ is homologous to $\beta$'' we mean that $T\beta - \beta$ is
homologous to zero as a closed $1$-chain.\\
\\
\begin{proof}
We will prove the lemma in the case where the intersection graph $G$ is a tree (i.e.
connected). The general case may be proved with a small variation.\\
Using lemma~\ref{lem:1}, for each $\pi$ there is some $f:\Sigma \rightarrow
\Sigma$
, which is a composition of twists, such that
$T=f^{-1} \circ T_\pi \circ f$ . If the lemma is
 true for $T_\pi$, and $\beta$ is as above, then $T\beta = (f^{-1}\circ
T_\pi \circ f) \beta$. From the assumption, one verifies that
$T_\pi [\theta]
=[\theta]$ for $\theta = f \beta$. Then the lemma applied to $T_\pi$ implies
that $f[\beta]$ has intersection number zero with all $[\alpha _{i}]$s. Since
$f$ is a combination of twists along $\alpha_i$s, it follows that the lemma is
also true for $T$ . In other words, it is enough to verify the lemma for a
specific choice of the order of the twists.\\

To choose the appropriate order, note that the loops corresponding to
the vertices which are not leaves of $G$ form an independent subset of
$H_1(\Sigma ,\mathbb{Q})$ as a vector space over $\mathbb{Q}$. Extend this
subset to a basis of the subspace
generated by all $\alpha_i$s in $H_1(\Sigma,\mathbb{R})$ by inserting some
of the loops corresponding to
the leaves of $G$. Denote the elements of this basis by $\beta_1,...,\beta_k$
and the rest of $\alpha_i$s by $\gamma_1,...,\gamma_l$.
 Let $\eta_1,...,\eta_r$ be the set of curves among $\beta_j$s which are associated
with vertices of $G$ that are  neighbors of the leaves of $G$.
The remaining vertices (those that are not among $\gamma_i$s and $\eta_j$s)
may be divided  into three groups
$A,B,C$ as follows: $C$ is the group of leaves of $G$ among $\{\beta_1,...,\beta_k \}$.
The remaining vertices are divided into two sets $A,B$ such that
no element of $A$ is adjacent to an element of
$B$. We will first do the twists along elements of $A$, then along elements of
$B$ and then elements of $C$. Afterwards the twists along $\eta_j$s ,
which we will call the group $D$, are done,
and finally we do the twists along the elements of group $E$ which are
$\gamma_i$s. Let $\delta_1,...,\delta_n$
be the same as $\beta_1,...,\beta_k,\gamma_1,...,\gamma_l$ but with the new
ordering and accordingly,
 suppose that $\Gamma=[\gamma_1,...,\gamma_l]=[\delta_1,...,\delta_k]A^T=\Delta A^T$.\\

First, suppose that $\beta$ is a closed loop.
If $\lambda_i=(\beta |\delta_i), \theta_j=(\beta|\gamma_j)$, then
$\Theta = A \Lambda$,
and one should show that $\Lambda=0$, where the $i$th entry of $\Lambda$ is
$\lambda_i$ and the $j$th entry of $\Theta$ is $\theta_j$. \\

Consider the matrix $\Upsilon=[\epsilon_{ij}]$ with
$\epsilon_{ij}=1$ if $\delta_{i},\delta_{j}$ meet in a point and $i<j$ (otherwise,
$\epsilon_{ij}$ is defined to be zero). Since we are considering positive Dehn
twists, in the level of homology, $S_i$ (the twist along $\delta_i$) is
given by:
$S_i[\delta]=[\delta]+(\delta |\delta_i).[\delta_i]$.\\
It is then easy to check that:
\begin{equation}
[\beta]=S[\beta]=[\beta]+
\left(\begin{array}{cc}  \Lambda^{T} &  \Theta^{T}
\end{array}\right)
(I+\Upsilon+\Upsilon^{2}+ ...)
\left(\begin{array}{c} \delta_1\\
                       \vdots  \\
                       \delta_n
\end{array}\right) .
\end{equation}

Replacing $\Theta = A \Lambda$, it is enough to show that
$B=\left(\begin{array}{cc} I & A^{T}
\end{array}\right)
(I-\Upsilon)^{-1}
\left(\begin{array}{c} I\\
                       A
\end{array}\right)$
is invertible to get $\Lambda, \Theta=0$ (since $[\delta_1],...,[\delta_k]$ are
independent).
The matrix
$(I - \Upsilon)$ has the form:
$\left(\begin{array}{cc} I-E & -X \\
                        O & I
\end{array}\right) $ where
$x_{ij}$ is $1$ iff $\delta_i$ is the unique vertex among $\eta_p$s that is adjacent
to $\gamma_i$, and is zero otherwise. So
\begin{equation}
(I-\Upsilon)^{-1}=
\left(\begin{array}{ccc} (I-E)^{-1}&Y\\
                         0&I
\end{array}\right),
\end{equation}
with $X=(I-E)Y$. Hence
$B=(A^{T}+Y)A+(I-E)^{-1}$
which is invertible iff its multiplication on the left by $(I-E)$ is. \\

Decompose all of the matrices into the block corresponding to the groups $A,B,C,D$.
Then:
\begin{equation}
\begin{split}
\Delta &=[\Delta_A \ \ \  \Delta_B \ \ \  \Delta_C \ \ \ \Delta_D]=[\delta_1,...,\delta_k] \\
A&=[A_A \ \ \ A_B \ \ \ A_C \ \ \ 0]\\
X^T &=[0 \ \ \ 0 \ \ \ 0 \ \ \ X_D^T]\\
\end{split}
\end{equation}
\begin{displaymath}
E=\left[ \begin{array}{cccc}
0  &  \Delta_A^T .\Delta_B &  0  & \Delta_A^T .\Delta_D\\
0  &          0            &  0  & \Delta_B^T .\Delta_D\\
0  &          0            &  0  & \Delta_C^T .\Delta_D\\
0  &          0            &  0  & \Delta_D^T .\Delta_D
\end{array}
\right].
\end{displaymath}

Notice that $(\Delta_A^T .\Delta_B)A_B$  is zero since:
\begin{equation}
0= \Delta_A^T .\Gamma= (\Delta_A^T .\Delta)A^T= (\Delta_A^T .\Delta_B)A_B,
\end{equation}
where the last equality follows from the fact that
$ \Delta_A^T .\Delta_A ,\Delta_A^T .\Delta_C, A_D$ are all zero.
Consequently:\\
\begin{equation}
EA^T=[ (\Delta_A^T .\Delta_B)A_B \ \ \ 0 \ \ \ 0 \ \ \ 0]^T=0.
\end{equation}
This means that $(I-E)B=I+A^TA+XA$ which is always invertible for any  $A,X$ of the
above block form.\\

Now suppose that $\beta$ is not closed. Then we will have:\\
\begin{equation}
\Theta=\beta.\Gamma^T=A(\beta.\Delta^T)+R^T=A\Lambda + R^T,
\end{equation}
which means that we will get:
\begin{equation}
\begin{split}
S[\beta]=[\beta]&+\Lambda^T
\left(\begin{array}{cc}  I &  A^T
\end{array}\right)
(I-\Upsilon)^{-1}
\left(
\begin{array}{c} I\\A
\end{array}\right)
\left(\begin{array}{c} \delta_1\\
                       \vdots  \\
                       \delta_k
\end{array}\right)\\
&+R\left(\begin{array}{c}
\gamma_1 \\ \vdots \\ \gamma_l
\end{array}\right).
\end{split}
\end{equation}

Putting $Q=\Lambda^T (I-E)^{-1}$, this equation reads as:
\begin{equation}
0=Q(I+A^TA+XA)\Delta ^T + R \Gamma ^T.
\end{equation}

Put $Q=\Lambda^T(I-E)^{-1}$. Then $Q(I+A^TA+XA)+RA=0$. If $Q=[Q_A,Q_B,Q_C,Q_D]=[Q_0,Q_D]$
and similarly $A=[A_0, 0], X^T=[0,X_D^T]$ then this equation reads as:
\begin{equation}
[Q_0\ \ \ Q_D](I+
\left[\begin{array}{cc}
A_0^TA_0  & 0 \\
X_DA_0    & 0
\end{array}\right])+[RA_0 \ \ \ 0]=[* \ \ \ Q_D]=0.
\end{equation}

So $Q_D=0$ which implies that $QX=0$. Combining this we the previous computation $(I-E)^{-1}A^TA=A^TA$
(which is the same as $EA^TA=0$), it follows that:
\begin{equation}
\begin{split}
0&=(\beta .\Delta)((I-E)^{-1}+A^TA)+RA\\
&=(\beta .\Delta)(I+A^TA(I-E))+RA(I-E).
\end{split}
\end{equation}

Again using the block representation of the matrix $E$ and the easy fact that
$A_A\Delta_A^T\Delta_B=0$, we see that:
\begin{equation}
\begin{split}
AE&=[0 \ \ \ \ A_A\Delta_A^T.\Delta_B \ \ \ 0 \ \ \ \ A \Delta^T .\Delta_D]\\
&=[0 \ \ \ 0 \ \ \ 0 \ \ \ \Gamma^T.\Delta_D]=\Gamma^T.\Delta .\\
\end{split}
\end{equation}

The result is that
\begin{equation}\label{eq:1}
\begin{split}
0&=(\beta.\Delta)+((\beta .\Delta)A^T+R)(A-\Gamma ^T .\Delta)\\
&=(\beta.\Delta)+(\beta. \Gamma^T)(A-\Gamma^T .\Delta) .\\
\end{split}
\end{equation}

In particular by looking at the first three blocks we get:
\begin{equation}
0=[\beta.\Delta_A \ \ \ \beta.\Delta_B \ \ \ \beta.\Delta_C]
\left[\begin{array}{c}
\Delta_A^T\\
\Delta_B^T\\
\Delta_C^T
\end{array}\right] +(\beta. \Gamma^T)\Gamma.
\end{equation}

This means that $\sum_{\alpha \in \mathcal{A}}(\beta | \alpha)\alpha =0$, where $\mathcal{A}$
is the set of all $\delta_i$s in $A,B,C$ or $E$.\\
\\
Let $\mathcal{B}$ be the set of elements $\alpha \in \mathcal{A}$ such that $(\beta | \alpha)$
is not zero. Then we claim that the elements of $\mathcal{B}$ are disjoint loops on $\Sigma$.
The reason is that the associated graph of the elements of  $\mathcal{B}$ is again a forest.
If it has any edges then we may find elements $\eta, \mu \in  \mathcal{B}$ such that they
are connected in the associated graph and $\eta$ is a leaf of this graph. Then:
\begin{equation}
0=\sum_{\alpha \in \mathcal{B}}(\beta | \alpha)(\alpha|\eta)=(\beta|\mu)(\mu|\eta)=\pm(\beta|\mu),
\end{equation}
which is a contradiction, since we had assumed that $(\beta | \alpha)\neq 0$
for all $\alpha \in \mathcal{B}$.\\
\\
\begin{figure}
\mbox{\vbox{\epsfbox{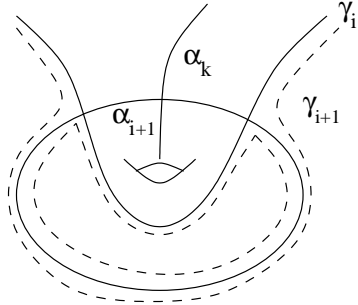}}}
\caption{\label{fig:3}
{$\gamma_{i+1}$ may be obtained from $\gamma_i$ by connecting couples of
intersection points with $\alpha_{i+1}$ along $\alpha_{i+1}$ and disjoint
 from $\alpha_k$ for $k \le i$}}
\end{figure}
\\
\begin{lem}\label{lem:4} Suppose that $\beta$ is a path on a surface $\Sigma$.
If $\hat{\gamma_1},...,\hat{\gamma_r}$ are disjoint loops
on $\Sigma$ with $(\beta|\hat{\gamma_i})\neq 0$ for all $i$ and
$\sum_{i}(\beta|\hat{\gamma_i})\hat{\gamma_i}=0$ in the homology,
then at least one of the $\hat{\gamma_i}$'s
is homologically trivial.\end{lem}
\begin{proof}(of lemma~\ref{lem:4}) Suppose that $\partial D=\sum_i(\beta|\hat{\gamma_i})
\hat{\gamma_i}$ and that $D=\sum_j n_j D_j$ with $D_j$s domains on $\Sigma$
and $n_i\neq n_j$ for $i \neq j$. Let $n_1>...>n_s$ and suppose that $n_1>0$
(changing $\beta$ with $-\beta$ if necessary, this may be assumed). Since
 each boundary component is common between at most two of the domains and
$n_1$ is the biggest of all $n_j$s, $\partial D_1=\sum_i a_i \hat{\gamma_i}$
(i.e. the boundary of $D_1$ is a combination of the loops) and $a_i$ has the
same sign as $(\beta|\hat{\gamma_i})$ .\\
Since $D_1$ is a domain, $(\beta|\partial D_1)\leq 1$. On the other hand,
$(\beta|\partial D_1)=\sum_i a_i(\beta|\hat{\gamma_i})$ and
$a_i (\beta|\hat{\gamma_i}) >0$ unless $a_i=0$. Thus, exactly one of the
$a_i$s is nonzero and one of the loops should be homologically
trivial.\end{proof}

Since we have assumed that non of $\alpha_i$s is homologically trivial, it should be the
case that $\mathcal{B}$ is empty, i.e. $(\beta|\alpha)=0$ for any $\alpha$ in the groups
$A,B,C,E$. In particular $\beta .\Gamma^T=0$ which implies by equation~\ref{eq:1} that
$\beta. \Delta^T=0$. This
completes the proof of lemma~\ref{lem:3}.\end{proof}
\begin{proof}(of proposition~\ref{prop:disks}) Put $f=T^{\epsilon}$
and suppose that there exists a path $\gamma$ from $x_-$ to $x_+$ such that
$[\gamma]-[f\gamma]$ represents the trivial homology class in the homology of
$\Sigma$. By lemma~\ref{lem:3},
$[\gamma]$ has zero intersection number with all $[\alpha_i]$s. \\
\\
\begin{lem}\label{lem:5} Suppose that $\gamma$ is a path from $x_-$ to $x_+$ such that
it has zero intersection number with all $\alpha_i$s. Then there is a
union $\delta=\cup_j \delta_j$
of paths $\delta_j \subset$Fix$(T)$
such that $[\delta]
-[\gamma]$ represents the trivial homology class in $\Sigma$.\end{lem}
\begin{proof}
After changing
indices if necessary, we may assume that $\alpha_i$ intersects at most
one $\alpha_k$ for $k<i$. Define $\gamma_i$ as follows: $\gamma_0=\gamma$.
Suppose that $\gamma_i$ is defined such that it does not cut $\alpha_1,...,\alpha_i$
and it has zero intersection number with $\alpha_{i+1},\alpha_{i+2},...$ .
One may divide the intersection points of $\alpha_{i+1}$ and $\gamma_i$
into couples of positive and negative intersections, say $p_{\pm}^1,...,p_{\pm}^k$.
Cut $\gamma_i$ at all $p_{\pm}^l$s to get several paths on $\Sigma$. Two of the end
points of these paths are near $p_+^l$ and two are near $p_-^l$. Since
$\alpha_{i+1}\setminus \cup_{j=1}^{i}\alpha_j$ is connected, we may connect these
two pairs by two other paths parallel to $\alpha_{i+1}\setminus \cup_{j=1}^{i}\alpha_j$.
If we do this for $l=1,...,k$ we will get a new $1$-chain $\gamma_{i+1}$ which does not
cut $\alpha_1,....,\alpha_{i+1}$
(see figure~\ref{fig:3}).
 $[\gamma_{i+1}]-[\gamma_i]$ is definitely homologically trivial,
hence so is $[\gamma_{i+1}]-[\gamma]$. $\delta=\gamma_n$ will then be
separated from all $\alpha_i$s and may be assumed to lie in Fix$(T)$.\end{proof}

Since $x_-,x_+$ will be the endpoints of one of the connected pieces of
$\delta=\cup_j \delta_j$, it is automatically implied that they are on
the same $V_i$. This finishes the proof of proposition~\ref{prop:disks}\end{proof}

\begin{proof}(of proposition~\ref{prop:energy})
For an element $u \in \mathcal{M}(\varphi)$, the energy $E(u)$ is defined by:
\begin{equation}
E(u)=\int_{\mathbb{R}\times [0,1]} |du|^2=
\int_{\mathbb{R}\times [0,1]} u^*\omega_R,
\end{equation}
where $\varphi \in \pi_{2}(x_-,x_+)$ and the metric is given using $\omega_R$ and $J$.\\
If $x_-,x_+$ are in the same piece and $\gamma(s)=u(s,1)$ then $T\gamma$ is
homologous to $\gamma$, which implies by proposition~\ref{prop:disks} and lemma~\ref{lem:5} that
there is a $1$-chain $\delta$ in $Fix(T)$ homologous to $\gamma$.
Define $v_j(s,t)=(H^{(1-t)\epsilon}(\delta_j(s)),t)$ where $\delta_j$s
are different pieces of $\delta$. If $L$ is the $2$-chain with
$\partial [L]=[\gamma]-[\delta]$, and $\Tilde{u}(s,t)=(u(s,t),t)$ then
\begin{displaymath}
[K]=[L\times \{1\}]-(\sum_{j}[v_j])+[\Tilde{u}]-[f(L)\times \{0\}]
\end{displaymath}
defines a homology class in $H_2(\Sigma \times \mathbb{R})=H_2(\Sigma)$. So
\begin{equation}
[K]=n(\varphi)[\Sigma].
\end{equation}

Note that $n(\varphi)$ does not depend on the representative $u$ of $\varphi$.
In fact if $u_1,u_2$ are representatives of the classes $\varphi_1,\varphi_2
\in {\pi_2}(x_-,x_+), u_i \in \mathcal{M}(\varphi_i)$ and
$[K_1]=n_1[\Sigma],[K_2]=n_2[\Sigma]$
are the associated classes, then  $n_1 - n_2=[\{z\}\times [0,1]].
([\Tilde{u_1}]-[\Tilde{u_2}])= [\{z\}\times [0,1]].
([\Tilde{\varphi_1}]-[\Tilde{\varphi_2}])$ and in particular if $\varphi_1=
\varphi_2$ then $n_1=n_2$ (here, $\Tilde{g}(s,t)=(g(s,t),t)$ for all $g$).\\

If we denote the pull back via
the first projection $p_1:\Sigma \times \mathbb{R}\rightarrow \Sigma$ of
$\omega_R$ by ${\omega}_R$ as well, then:
\begin{equation}
\begin{split}
\int_{[K]}{\omega}_R & =\int_{\mathbb{R}\times [0,1]}u^*\omega_R
+\int_{L}\omega_R -\int_{L}f^*\omega_R -
\sum_j \int_{\mathbb{R} \times [0,1]} {v_j}^* \omega_R \\
 & =E(u) - \sum_{j}\int_{\mathbb{R} \times [0,1]} \omega_R(dv_j\frac{\partial}
{\partial s}, \epsilon X_{h}(v_j))\\
 & = E(u)-\epsilon (h(x_+)-h(x_-)).\\
\end{split}
\end{equation}

Hence,
\begin{equation}
E(u)=n(\varphi)(\int_{\Sigma}
\omega_R)
+\epsilon (h(x_+)-h(x_-)),
\end{equation}
which completes the proof of the proposition.\end{proof}

\section{Main Theorem}
\label{sec:4}
The content of this section is the proof of the following computation:\\
\\
\begin{thm}\label{thm:main} Let $(\Sigma,\alpha_1,...,\alpha_n)$ be an acceptable
setting and put $C=\cup_i \alpha_i$. Let $\psi$ be the isotopy class of the
combination $T$ of the positive Dehn twists along $\alpha_i$s in some order in the
mapping class group  $\Gamma=\pi_0(Diff^+(\Sigma))$ .
 Then $HF^*(\psi)=H^*(\Sigma,C)$ as
$\mathbb{Z}/2$-graded
$H^*(\Sigma)$-modules, where $H^*(\Sigma)=HF^*(Id)$ acts on the right hand
side by the ordinary cup product and on the left hand side by quantum cup product.\end{thm}

 For a diffeomorphism
$f:\Sigma \rightarrow \Sigma$, let $T_f\Sigma$ be the mapping
torus of $f$ defined by
\begin{equation}
T_f\Sigma:=\frac{\Sigma \times I}{\sim} \ \ \ \ (x,0) \sim (f(x),1)
\ \ \ \ \forall x \in \Sigma.
\end{equation}

For $f\in Symp(\Sigma,\omega)$, the group of symplectomorphisms of $\Sigma$
with respect to the symplectic form $\omega$,
let $\Tilde{\omega}$ be the pull back of
$\omega$ to $\Sigma \times \mathbb{R}$ and denote the induced form on
$T_f(\Sigma)$ by $\Tilde{\omega}$ as well. $T_f(\Sigma)$ is fibered over
$S^1$. The Euler class of the tangent plane bundle along the fibers of
$T_f(\Sigma)$, will be an element $c_f$ of $H^2(T_f(\Sigma),\mathbb{R})$.
Consider the exact sequence:\\
\begin{equation}
...\rightarrow H^1(\Sigma,\mathbb{R}) \xrightarrow{Id-f^*}
H^1(\Sigma,\mathbb{R}) \xrightarrow{d} H^2(T_f(\Sigma),\mathbb{R})
\rightarrow H^2(\Sigma,\mathbb{R}) \rightarrow 0.
\end{equation}

Since $[\Tilde{w}]-\frac{\int_{\Sigma}\omega}{\chi(\Sigma)}[c_f]$ vanishes
when restricted to the fibers, it will be of the form $d(m(f))$ with
$m(f) \in H^1(\Sigma,\mathbb{R})/Im(Id-f^*)$. $m(\varphi)$ is called
the \emph{monotonicity class} and
$f$ is called \emph{monotone} if
$m(f)=0$ or saying it differently, if $[\Tilde{\omega}]=\frac{\int_{\Sigma}
\omega}{\chi(\Sigma)}[c_f]$.
It is proved in \cite{Sei2} that if two monotone symplectomorphisms represent
the same isotopy class $[\psi]$ in the mapping class group $\Gamma$ then
they will give the same Floer cohomology denoted by $HF^*([\psi])$. It is also
shown that there are monotone representatives in each isotopy class.\\
\\
\begin{lem}\label{lem:6} $T^{\epsilon}$ is monotone.\end{lem}
\begin{proof} Put $f=T^{\epsilon}$. We should show that for a class
$[\gamma] \in H_1(\Sigma,\mathbb{R})$ with $[f_*(\gamma)]=[\gamma]$,
$\eta=\Tilde{\omega}-\frac{\int_{\Sigma}
\omega}{\chi(\Sigma)}[c_f]$ evaluates zero on the homology class associated to
$[\gamma]$ in $H_2(T_f(\Sigma),\mathbb{R})$.\\
If for some two chain $L_0$, $\partial [L_0]=[\gamma]-[f(\gamma)]$ then the
class in $H_2(T_f(\Sigma))$ associated to $[\gamma]$ is the image of
$[L_1]=[\gamma \times [0,1]]-[L_0\times \{0\}]$ under the quotient map
$\pi:\Sigma \times [0,1]\rightarrow T_f(\Sigma)$.\\

By lemma~\ref{lem:3} $\gamma$ has zero intersection with $\alpha_i$s and
by lemma~\ref{lem:5} (applied to $x_-=x_+$ on $\gamma$) it is homologous
to a sum $\sum_{j}\delta_j$ of loops
fixed by $T^R$. Let $\partial [L_2]=[\gamma]-\sum_{j}[\delta_j]$ and define
$v_j$s as in the proof of proposition~\ref{prop:energy}. Then:
\begin{equation}
[K]=[L_2\times \{1\}]-(\sum_{j}[v_j])+L_1-[f(L_2)\times \{0\}]
\end{equation}
will be a homology class in $H_2(\Sigma \times \mathbb{R})=H_2(\Sigma)$.
The pull back of $\eta$ evaluates to zero
on this class. Since $\omega_R$ is invariant under $f$, and $c_f$ comes from
the cohomology of $T_f(\Sigma)$, the evaluation of
$\eta$ on $[L_2 \times \{1\}]-[f(L_2)\times \{0\}]$ will be zero as well.
One can choose a trivialization of the $2$-plane bundle in a neighborhood of
$\delta_j$ and pull this trivialization back via $H^{(t-1)\epsilon}$ to
$v_j(.,t)$ to get a trivialization along $[v_j]$. So
$c_f$ is trivial on $[v_j]$s.  Finally the same argument as that of
proposition~\ref{prop:energy} implies that the integration of
$\Tilde{\omega}^R$ on the $2$-chain $[v_j]$ is zero.\\
Consequently $\eta$ is zero on $[L_1]\in H_2(T_f(\Sigma),\mathbb{R})$.\end{proof}

\begin{proof}(of theorem~\ref{thm:main})
$\pi_2(\Sigma)$ is trivial, hence there is no bubbling off
of holomorphic spheres. So the usual problem of non-compactness of moduli
space is automatically solved in our case.\\

Since $T^{\epsilon}$ is a monotone representative of $[\psi]$ it can be
used for the computation.
By proposition~\ref{prop:disks} there are no disks between the fixed points
in different pieces. Suppose that $x_-,x_+$ are two fixed points in the same
piece $V_j$.  We will show that all
elements $u$ of $\mathcal{M}(\varphi)$ (where $\varphi \in
\pi_2(x_-,x_+)$, $\mu(\varphi)=1$) satisfy $h(u(s,t))<3$ for all
$s,t$.\\

If we connect $x_-,x_+$ by a path $\gamma$ and put $\varphi_0(s,t)
=H^{(1-t)\epsilon}(\gamma(s))$, then $E(\varphi_0)=\epsilon(h(x_+)-h(x_-))$
by proposition~\ref{prop:energy}. On the other hand, from \cite{Sei2} we know
that since $T^{\epsilon,R}$ is monotone, $\mu(\varphi)-\mu(\varphi_0)=
\frac{\chi(\Sigma)}{\int_{\Sigma}\omega_R}(E(\varphi)-E(\varphi_0))$,
since $T^{\epsilon,R}$ is monotone. Hence
$\mu(\varphi)=\mu(\varphi_0)+n(\varphi)\chi(\Sigma)$.
Since $\mu(\varphi_0)$ is the Morse-index difference of $x_-$ and $x_+$ and
$\chi (\Sigma)<0$,
the only possibilities for classes $\varphi$ which contribute to equation (4)
satisfy $n(\varphi)=0$.\\

We claim that if $u \in \mathcal{M}(\varphi)$ for such $\varphi$, then the
image of $u$ should stay in $h^{-1}(-\infty,3]$. Suppose that this is not true.
Then we may find $s_0<s_1 , t_0<t_1$ such that for $U=[s_0,s_1]\times [t_0,t_1]
\subset \mathbb{R}\times [0,1]$ and any $(s,t)\in U$, $h(u(s,t))\in
I_0=[3-\delta,3+\delta]$. Then:
\begin{equation}
\begin{split}
\epsilon(h(x_+)-h(x_-))=\int_{\mathbb{R}\times [0,1]} |du|^2 &\geq \int_{U} |du|^2\\
                           &=\int_{U} u^*\omega_R=R\int_U u^*\omega,\\
\end{split}
\end{equation}
which is a contradiction since we may assume that $R$ is arbitrarily large.
This proves the claim.\\

Put $S=h^{-1}(-\infty,3]$ and note that the Hamiltonian flow maps $S$ to itself.
The Floer homology
of this  flow will give us the Floer homology of the map $T$ by the above discussion.
Then [8, theorem 7.3] shows that
$\mathcal{M}(x_-,x_+)$
is homeomorphic to the space of $u:\mathbb{R} \rightarrow S$ such that
\begin{equation}
\frac{du}{dt} = - \nabla _J h(u) \ \ \ \ with  \ \ \lim_{t \rightarrow \pm
\infty} u(t) = x_{\pm}.
\end{equation}
for a generic choice of $J$, this gives the isomorphism we are looking for.\end{proof}

\begin{figure}
\mbox{\vbox{\epsfbox{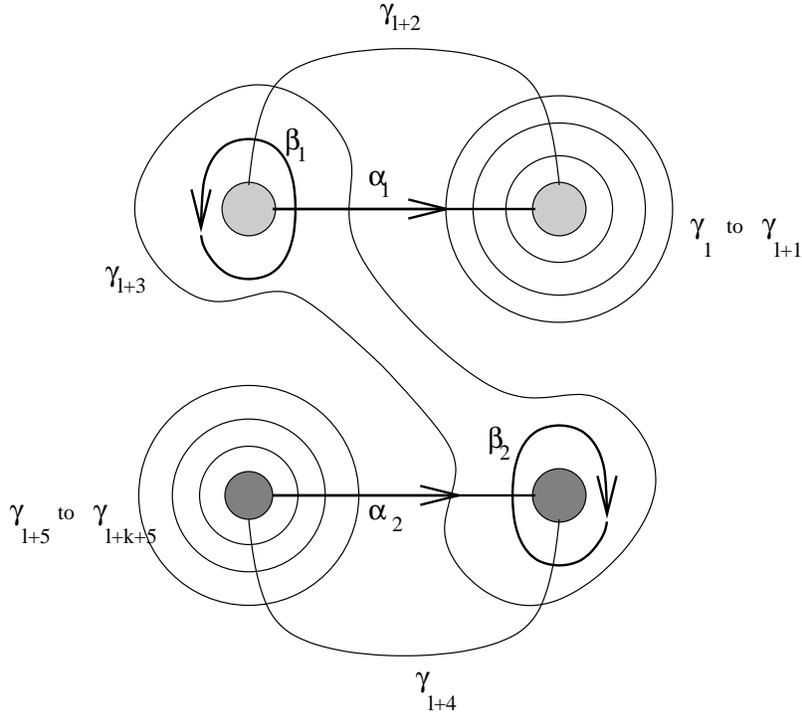}}}
\caption{\label{fig:4}
{The map $T$ is the combination of positive twists along $\gamma_i$s }}
\end{figure}
\begin{example}\label{ex:dehn.1}\end{example} On a surface of genus $2$, suppose that the twists
are done along the curves $\gamma_i$  as shown in
figure~\ref{fig:4}. Here the two pairs of shaded circles denote the attaching
circles of the two $1$-handles.
If $T=T_{l+k+5}\circ T_{l+k+4}\circ ... \circ T_1$ where $T_i$ is the
positive twist along $\gamma_i$, then the matrix of the action of $T$ on the
first homology of $\Sigma$ in the basis $\{[\alpha_1],[\beta_1],[\alpha_2],
[\beta_2]\}$ will be:
\begin{equation}
\left( \begin{array}{cccc}
0&0&-1&1\\
-k&-k&kl&-1\\
1&0&-l&0\\
1&1&-l&0\\
\end{array}\right).
\end{equation}

The characteristic polynomial will be $\chi_T(\zeta)=\zeta^2(\zeta+k)(\zeta+l)
+(\zeta+k)(\zeta+l)+(1-kl)$. This polynomial is irreducible over $\mathbb{Z}$ if
$k\equiv l$ $($mod $4)$ and they are both odd (by ``Eisenstein Criterion''
applied to $\chi_T(\zeta+1), p=2$).
If $\chi_T(\zeta)$ vanishes at some root of unity $\eta$, then it will be
divisible by $\zeta -(\eta+\overline{\eta})\zeta+1$ and $\eta+\overline{\eta}$
is a real number of absolute value at most $2$. If $k,l>4$ , it can be
easily seen that this is not the case. We quote the following two theorems from
\cite{CB}:\\
\\
\textbf{Theorem A.} \emph{If for $f:\Sigma \rightarrow \Sigma$, the characteristic
polynomial $\chi_f(\zeta)$ of the action of $f$ on the first homology, is
irreducible over $\mathbb{Z}$, has no roots of unity as zeros, and is not a
polynomial in $\zeta^n$ for any $n>1$, then $f$ is irreducible and
non-periodic.}
\\
\\
\textbf{Theorem B.} \emph{Every non-periodic irreducible automorphism of a closed
orientable hyperbolic surface is isotopic to a pseudo-Anosov automorphism.}
\\

Accordingly, the map $T$ will be representing a
pseudo-Anosov class and theorem~\ref{thm:main} explicitly computes the Floer homology
associated to its class in the mapping class group.
\\
\\
\begin{remark} Note that lemma~\ref{lem:3} remains true if $T$ is a
combination of positive twists along $\alpha_i$s and negative twists along
$\alpha_j'$s where the associated graphs $G,G'$ of $\alpha_i$s and
$\alpha_j'$s are forests, no $\alpha_i$ meets an $\alpha_j'$ and none of the
$\alpha_i$s is homotopic to a $\alpha_j'$.\end{remark}

\begin{defn} $(\Sigma, (\alpha_1,...,\alpha_n),(\alpha_1',...,
\alpha_m'))$ is called a \emph{strongly acceptable} setting if
$\alpha_i$s and $\alpha_j'$s are as above and non of $\alpha_i$s or
$\alpha_j'$s is homologically trivial.\end{defn}

The procedure of proposition~\ref{prop:morse} may be applied ``upside down'' to
give a Morse function $h$ with $h(\alpha_i)\subset (4,\infty)$ and
$h(\alpha_j') \subset (-\infty,-4)$ with the similar properties. For $T$ as
above, $T^\epsilon$ may be defined as the composition of $T$ with the
Hamiltonian flow of $h$. Its fixed points will be the critical points $p$
of $h$ with $|h(p)|<2$. Since the rest of our results do not remember the
sign of the twists, one may prove the
following theorem similar to theorem~\ref{thm:main} (c.f. the result of
\cite{Sei1}):\\
\\
\begin{thm}\label{thm:main.general} If $(\Sigma, (\alpha_1,...,\alpha_n),(\alpha_1',...,
\alpha_m'))$  is a strongly acceptable setting and $T$ is the composition
of positive Dehn twists along $\alpha_i$s and negative twists along
$\alpha_j'$s, in some order, then $HF^*(T)\simeq H^*(\Sigma \setminus
(\cup_i \alpha_i),(\cup_j \alpha_j'))$.\end{thm} $\square$

\begin{remark} Consider the quantum cup product\\
\begin{displaymath}
*: H^*(\Sigma,\mathbb{Z}/2) \otimes HF^*(T) \ \ \longrightarrow \ \ HF^*(T).
\end{displaymath}
The theorem just proven states that this map is the same as the usual
cup product\\
\begin{displaymath}
\cup : H^*(\Sigma, \mathbb{Z}/2) \otimes H^*(\Sigma,C) \ \ \longrightarrow
\ \ HF^*(\Sigma,C).
\end{displaymath}
Clearly, $H^2(\Sigma,\mathbb{Z}/2)$ acts trivially unless $C$ is empty (i.e.
$T$ is the identity class) and for $a \in H^1(\Sigma,\mathbb{Z}/2)$, if the
action on $H^*(\Sigma,C)$ is nonzero, then $a$ is dual to a curve $\mathit{l}
: S^1\rightarrow \Sigma$ which does not cut $C$. This is to say
that $T(\mathit{l})\simeq \mathit{l}$. These are the special cases of
{[10, theorem 1,2]} which are immediate because of the above computation.\end{remark}

\end{document}

%% file: command.tex
\newcommand{\divis}{\mathfrak d}
\newcommand\Dehn{D}
\newcommand\PiRed{\Pi^\red}

\newcommand\liftalpha{\widetilde\alpha}
\newcommand\liftbeta{\widetilde\beta}
\newcommand\sj{\mathfrak j}
\newcommand\MT{t}
\newcommand\Jacobi[2]{\left(\begin{array}{c} #1 \\ #2 \end{array}\right)}
\newcommand\Cobord{{\mathfrak C}(3)}
\newcommand\KOver{K_+}
\newcommand\KUnder{K_-}
\newcommand\Fiber{F}
\newcommand\prefspinc{\mathfrak u}
\newcommand\spincu{\mathfrak u}
\newcommand\SpinCCobord{\theta^c}
\newcommand\Field{\mathbb F}
\newcommand\Fmix[1]{F^{\mix}_{#1}}
\newcommand\ugc[1]{{\underline g}^\circ_{#1}}
\newcommand\uGc[1]{{\underline G}^\circ_{#1}}
\newcommand\uginf[1]{{\underline g}^\infty_{#1}}
\newcommand\uGinf[1]{{\underline G}^\infty_{#1}}
\newcommand\uec[1]{{\underline e}^\circ_{#1}}
\newcommand\uEc[1]{{\underline E}^\circ_{#1}}
\newcommand\Ec[1]{{E}^\circ_{#1}}
\newcommand\Gc[1]{{G}^\circ_{#1}}
\newcommand\ueinf[1]{{\underline e}^\infty_{#1}}
\newcommand\uEinf[1]{{\underline E}^\infty_{#1}}
\newcommand\rspinc{\underline\spinc}
\newcommand\Fc{F^\circ}
\newcommand\fc{f^\circ}
\newcommand\uFc{{\underline F}^\circ}
\newcommand\ufc{{\underline f}^\circ}
\newcommand\uF{{\underline F}^\circ}
\newcommand\uHFc{\uHF^\circ}
\newcommand\uCFc{\uCF^\circ}
\newcommand\HFc{\HF^\circ}
\newcommand\CFc{\CF^\circ}

\newcommand\Dual{\mathcal D}
\newcommand\Duality\Dual
\newcommand\aPrBd{\partail_{\widetilde \alpha}}
\newcommand\bPrBd{\partail_{\widetilde \beta}}
\newcommand\cPrBd{\partail_{\widetilde \gamma}}

\newcommand\liftTa{{\widetilde{\mathbb T}}_\alpha}
\newcommand\liftTb{{\widetilde{\mathbb T}}_\beta}
\newcommand\liftTc{{\widetilde{\mathbb T}}_\gamma}
\newcommand\liftTd{{\widetilde{\mathbb T}}_\delta}
\newcommand\liftTcPr{{\widetilde{\mathbb T}}_\gamma'}

\newcommand\Tor{\mathrm{Tor}}
\newcommand\TaPr{\Ta'}
\newcommand\TbPr{\Tb'}
\newcommand\TcPr{\Tc'}
\newcommand\TdPr{\Td'}
\newcommand\alphaprs{\alphas'}
\newcommand\betaprs{\betas'}
\newcommand\gammaprs{\gammas'}
\newcommand\Spider{\sigma}
\newcommand\EulerMeasure{\widehat\chi}

\newcommand\RelInv{F}

\newcommand\Knot{\mathbb K}
\newcommand\Link{\mathbb L}
\newcommand\SubLink{\mathbb K}
\newcommand\SublinkC{{\SubLink}'}
\newcommand\liftSym{{\widetilde \Sym}^g(\Sigma)_\xi}
\newcommand\liftCFinf{\widetilde{CF}^\infty}
\newcommand\liftz{\widetilde z}
\newcommand\liftCFp{\widetilde{CF}^+}
\newcommand\liftCF{\widetilde{CF}}
\newcommand\liftHFinf{\widetilde{HF}^\infty}
\newcommand\liftHFp{\widetilde{HF}^+}
\newcommand\liftx{\widetilde\x}
\newcommand\lifty{\widetilde\y}
\newcommand\liftSigma{\widetilde\Sigma}
\newcommand\liftalphas{\widetilde \alphas}
\newcommand\liftbetas{\widetilde \betas}
\newcommand\liftgammas{{\widetilde\gammas}}
\newcommand\liftdeltas{{\widetilde\deltas}}
\newcommand\liftgammasPr{\liftgammas'}
\newcommand\liftGamma{\widetilde{\Gamma}}
\newcommand\liftdeltasPr{\liftdeltas'}
\newcommand\Tz{{\widetilde{\mathbb T}}_\zeta}
\newcommand\liftTheta{\widetilde\Theta}

\newcommand\liftF[1]{{\widetilde F}_{#1}}
\newcommand\spinck{\mathfrak k}
\newcommand\Laurent{\mathbb L}
\newcommand\renEuler{\widehat\chi}
\newcommand\spinccanf{k}
\newcommand\spinccan{\ell}

\newcommand\mCP{{\overline{\mathbb{CP}}}^2}

\newcommand\HFpRed{\HFp_{\red}}
\newcommand\HFmRed{\HFm_{\red}}
\newcommand\mix{\mathrm{mix}}

\newcommand\RP[1]{{\mathbb{RP}}^{#1}}

\newcommand\liftGr{\widetilde\gr}
\newcommand\chiTrunc{\chi^{\mathrm{trunc}}}
\newcommand\chiRen{\widehat\chi}
\newcommand\Fred[1]{F^{\mathrm{red}}_{#1}}

\newcommand{\F}{\mathbb F}
\newcommand\SpinCz{\mathfrak S}
\newcommand\interior{\mathrm{int}}

\newcommand{\ev}{\mathrm{ev}}
\newcommand{\odd}{\mathrm{odd}}
\newcommand\sRelSpinC{\underline\spinc}
\newcommand\RelSpinC{\underline{\SpinC}}
\newcommand\relspinc{\underline{\spinc}}
\newcommand\ThreeCurveComp{\Sigma-\alpha_{1}-\ldots-\alpha_{g}-\beta_{1}-\ldots-\beta_{g}-\gamma_1-...-\gamma_g}
\newcommand\SpanA{{\mathrm{Span}}([\alpha_i]_{i=1}^g)}
\newcommand\SpanB{{\mathrm{Span}}([\beta_i]_{i=1}^g)}
\newcommand\SpanC{{\mathrm{Span}}([\gamma_i]_{i=1}^g)}
\newcommand\Filt{\mathcal F}
\newcommand\HFinfty{\HFinf}
\newcommand\CFinfty{\CFinf}
\newcommand\Tai{{\mathbb T}_{\alpha}^i}
\newcommand\Tbj{{\mathbb T}_{\beta}^j}
\newcommand\RightFp{R^+}
\newcommand\LeftFp{L^+}
\newcommand\RightFinf{R^\infty}
\newcommand\LeftFinf{L^\infty}
\newcommand\Area{\mathcal A}
\newcommand\PhiIn{\phi^{\mathrm{in}}}
\newcommand\PhiOut{\phi^{\mathrm{out}}}

\newcommand\x{\mathbf x}
\newcommand\w{\mathbf w}
\newcommand\z{\mathbf z}
\newcommand\p{\mathbf p}
\newcommand\q{\mathbf q}
\newcommand\y{\mathbf y}
\newcommand\Harm{\mathcal H}
\newcommand\sS{\mathcal S}
\newcommand\sW{\mathcal W}
\newcommand\sX{\mathcal X}
\newcommand\sY{\mathcal Y}
\newcommand\sZ{\mathcal Z}
\newcommand\cent{\mathrm{cent}}

\newcommand\BigO{O}
\newcommand\ModFlowMod{\ModFlow^{*}}
\newcommand\ModCent{\ModFlow^{\mathrm{cent}}
\left({\mathbb S}\longrightarrow \Sym^{g-1}(\Sigma_{1})\times 
\Sym^2(\Sigma_{2})\right)}
\newcommand\ModSphere{\ModFlow\left({\mathbb S}\longrightarrow
\Sym^{g-1}(\Sigma_{1})\times \Sym^2(\Sigma_{2})\right)}
\newcommand\ModSpheres\ModSphere
\newcommand\CF{CF}
\newcommand\cyl{\mathrm{cyl}}
\newcommand\CFa{\widehat{CF}}
\newcommand\CFp{\CFb}
\newcommand\CFm{\CF^-}
\newcommand\CFleq{\CF^{\leq 0}}
\newcommand\HFleq{\HF^{\leq 0}}
\newcommand\CFmeq{\CFleq}
\newcommand\CFme{\CFleq}
\newcommand\HFme{\HFleq}
\newcommand\HFpred{\HFp_{\rm red}}
\newcommand\HFpEv{\HFp_{\mathrm{ev}}}
\newcommand\HFpOdd{\HFp_{\mathrm{odd}}}
\newcommand\HFmred{\HFm_{\rm red}}
\newcommand\HFred{\HF_{\rm red}}
\newcommand\coHFm{\HF_-}
\newcommand\coHFp{\HF_+}
\newcommand\coHFinf{\HF_\infty}
\newcommand{\red}{\mathrm{red}}
\newcommand\ZFa{\widehat{ZF}}
\newcommand\ZFp{ZF^+}
\newcommand\BFp{BF^+}
\newcommand\BFinf{BF^\infty}
\newcommand\HFp{\HFb}
\newcommand\HFpm{HF^{\pm}}
\newcommand\HFm{\HF^-}
\newcommand\CFinf{CF^\infty}
\newcommand\HFinf{HF^\infty}
\newcommand\CFb{CF^+}
\newcommand\HFa{\widehat{HF}}
\newcommand\HFb{HF^+}
\newcommand\gr{\mathrm{gr}}
\newcommand\Mas{\mu}
\newcommand\UnparModSp{\widehat \ModSp}
\newcommand\UnparModFlow\UnparModSp
\newcommand\Mod\ModSp
\newcommand{\cpl}{{\mathcal C}^+}
\newcommand{\cmi}{{\mathcal C}^-}
\newcommand{\cplm}{{\mathcal C}^\pm}
\newcommand{\cald}{{\mathcal D}}

\newcommand\PD{\mathrm{PD}}
\newcommand\dist{\mathrm{dist}}
\newcommand\Paths{\mathcal{B}}
\newcommand\Met{\mathfrak{Met}}
\newcommand\Lk{\mathrm{Lk}}
\newcommand\Jac{\mathfrak{J}}
\newcommand\spin{\mathfrak s}

\newcommand{\xTuple}{\underline x}
\newcommand{\yTuple}{\underline y}
\newcommand{\aTuple}{\underline a}
\newcommand{\bTuple}{\underline b}
\newcommand{\xiTuple}{\underline \xi}
\newcommand{\spinc}{\mathfrak s}
\newcommand{\spincfour}{\mathfrak r}
\newcommand{\spincX}{\spincfour}
\newcommand{\spinct}{\mathfrak t}
\newcommand{\Exp}{\mathrm{Exp}}
\newcommand{\etaTuple}{\underline \eta}

\newcommand\Real{\mathrm Re}
\newcommand\brD{F}
\newcommand\brDisk{\brD}
\newcommand\brPhi{\widetilde\phi}
\newcommand\Perm[1]{S_{#1}}
\newcommand\BranchData{\mathrm Br}
\newcommand\ModMaps{\mathcal M}
\newcommand\ModSp\ModMaps
\newcommand\ModMapsComp{\overline\ModMaps}
\newcommand\tsModMaps{\ModMaps^\circ}
\newcommand\ModCurves{\mathfrak M}
\newcommand\ModCurvesComp{\overline\ModCurves}
\newcommand\ModCurvesSmoothInt{\ModCurves_{si}}
\newcommand\ModCurvesSmoothBoundary{\ModCurves_{sb}}
\newcommand\BranchMap{\Phi}
\newcommand\SigmaMap{\Psi}
\newcommand\ProdSig[1]{\Sigma^{\times{#1}}}
\newcommand\ProductForm{\omega_0}
\newcommand\ProdForm{\ProductForm}
\newcommand\Energy{E}

\newcommand\qTup{\mathfrak q}
\newcommand\pTup{\mathfrak p}

\newcommand\Cinfty{C^{\infty}}
\newcommand\Ta{{\mathbb T}_{\alpha}}
\newcommand\Tb{{\mathbb T}_{\beta}}
\newcommand\Tc{{\mathbb T}_{\gamma}}
\newcommand\Td{{\mathbb T}_{\delta}}
\newcommand\Na{N_{\alpha}}
\newcommand\Nb{N_{\beta}}
\newcommand\Torus{\mathbb{T}}
\newcommand\Diff{\mathrm{Diff}}

\newcommand\dbar{\overline\partial}
\newcommand\Map{\mathrm{Map}}
\newcommand\Strip{\mathbb{D}}

\newcommand\mpa{p}
\newcommand\mpb{q}
\newcommand\Mac{M}

\newcommand\del{\partial}

\newcommand\alphas{\mbox{\boldmath$\alpha$}}
\newcommand\mus{\mbox{\boldmath$\mu$}}
\newcommand\xis{\mbox{\boldmath$\xi$}}
\newcommand\etas{\mathbf\eta}
\newcommand\betas{\mbox{\boldmath$\beta$}}
\newcommand\gammas{\mbox{\boldmath$\gamma$}}
\newcommand\deltas{\mbox{\boldmath$\delta$}}
\newcommand\HFto{HF _{{\rm to}} ^{SW}}
\newcommand\CFto{CF _{{\rm to}} ^{SW}}
\newcommand\HFfrom{HF_{{\rm from}} ^{SW}}
\newcommand\CFfrom{CF_{{\rm from}} ^{SW}}

\newcommand\NumPts{n}
\newcommand\Ann{\mathrm{Ann}}

\newcommand\Dom{\mathcal D}
\newcommand\PerClass[1]{{\PerDom}_{#1}}
\newcommand\PerClasses[1]{{\Pi}_{#1}}
\newcommand\PerDom{\mathcal P}
\newcommand\RenPerDom{\mathcal Q}
\newcommand\intPerDom{\mathrm{int}\PerDom}
\newcommand\csum{*}
\newcommand\CurveComp{\Sigma-\alpha_{1}-\ldots-\alpha_{g}-\beta_{1}-\ldots-\beta_{g}}
\newcommand\EmbSurf{Z}
\newcommand\Mult{m}

\newcommand\NumDoms{m}

\newcommand\uCF{\underline\CF}
\newcommand\uCFinf{\uCF^\infty}
\newcommand\uHF{\underline{\HF}}
\newcommand\uHFp{\underline{\HF}^+}
\newcommand\uHFmred{\underline{\HF}^-_{\red}}
\newcommand\uHFpRed{\uHFpred}
\newcommand\uHFpred{\underline{\HF}^+_{\red}}
\newcommand\uCFp{\underline{\CFp}}
\newcommand\uHFm{\underline{\HF}^-}
\newcommand\uHFa{\underline{\HFa}}
\newcommand\uCFa{\underline{\CFa}}
\newcommand\uDel{\underline{\partial}}
\newcommand\uHFinf{\uHF^\infty}

\newcommand\alphaperps{\{\alpha_1^\perp,...,\alpha_g^\perp\}}
\newcommand\betaperps{\{\beta_1^\perp,...,\beta_g^\perp\}}
\newcommand\gammaperps{\{\gamma_1^\perp,...,\gamma_g^\perp\}}
\newcommand\ufinf{\underline{f^\infty_{\alpha,\beta,\gamma}}}
\newcommand\uHFleq{\underline{\HFleq}}
\newcommand\uCFleq{\underline{\CFleq}}
\newcommand\ufp[1]{\underline{f^+_{#1}}}
\newcommand\ufleq[1]{\underline{f^{\leq 0}_{#1}}}
\newcommand\ufa{\underline{{\widehat f}_{\alpha,\beta,\gamma}}}
\newcommand\uFa[1]{\underline{{\widehat F}_{#1}}}
\newcommand\uFm[1]{{\underline{F}}^-_{#1}}
\newcommand\Fstar[1]{F^{*}_{#1}}
\newcommand\uFp[1]{{\underline{F^+}}_{#1}}
\newcommand\uFinf[1]{{\underline F}^\infty_{#1}}
\newcommand\uFleq[1]{{\underline F}^{\leq 0}_{#1}}
\newcommand\uFstar[1]{{\underline F}^{*}_{#1}}
\newcommand\Fm[1]{F^{-}_{#1}}
\newcommand\FmRed[1]{F^{-}_{#1}}
\newcommand\Fp[1]{F^{+}_{#1}}
\newcommand\Fpm[1]{F^{\pm}_{#1}}
\newcommand\Fa[1]{{\widehat F}_{#1}}
\newcommand\Finf[1]{F^{\infty}_{#1}}
\newcommand\Fleq[1]{F^{\leq 0}_{#1}}
\newcommand\Tx{{\mathbb T}_\xi}
\newcommand\Ty{{\mathbb T}_\eta}
\newcommand\Ring{\mathbb A}
\newcommand\HFstar{\HF^*}
\newcommand\uLeftFinf{{\underline L}^\infty}
\newcommand\uRightFinf{{\underline R}^\infty}
\newcommand\orient{\mathfrak o}
\newcommand\phiCov{\widetilde \phi}
\newcommand\xCov{\widetilde x}
\newcommand\wCov{\widetilde w}
\newcommand\yCov{\widetilde y}
\newcommand\qCov{\widetilde q}
\newcommand\alphaCov{\widetilde \alpha}
\newcommand\betaCov{\widetilde \beta}
\newcommand\gammaCov{\widetilde \beta}
\newcommand\piCov{{\widetilde \pi}_2}

\newcommand\CFKlgeq{CFK^{<,\geq}}
\newcommand\spincrel\relspinc
\newcommand\relspinct{\underline{\mathfrak t}}

\newcommand\CFKsa{CFK^0}
\newcommand\bCFKsa{^bCFK^0}
\newcommand\spincf{\mathfrak r}

\newcommand\Mass{\overline n}
\newcommand\Vertices{\mathrm{Vert}}

\newcommand\CFB{CFB}
\newcommand\CFBinf{\CFB^\infty}
\newcommand\Bouq{\mathbb B}
\newcommand\zees{\mathbf z}

\newcommand\MCone{M}
\newcommand\bCFKp{{}^b\CFKp}
\newcommand\bCFKm{{}^b\CFKm}
\newcommand\bPhip{{}^b\Phi^+}
\newcommand\bPhim{{}^b\Phi^-}
\newcommand\Phip{\Phi^+}
\newcommand\Phim{\Phi^-}
\newcommand\bCFKa{{}^b\CFKa}
\newcommand\CFK{CFK}
\newcommand\HFK{HFK}

\newcommand\CFKc{\CFK^\circ}
\newcommand\CFKp{\CFK^+}
\newcommand\CFKa{\widehat\CFK}
\newcommand\CFKm{\CFK^-}
\newcommand\CFKinf{\CFK^{\infty}}
\newcommand\HFKp{\HFK^+}
\newcommand\bHFKp{{}^b\HFK^+}
\newcommand\bHFKm{{}^b\HFK^-}
\newcommand\HFKa{\widehat\HFK}
\newcommand\HFKm{\HFK^-}
\newcommand\HFKinf{\HFK^{\infty}}
\newcommand\HFKc{\HFK^\circ}

\newcommand\Mark{m}
\newcommand\Marks{\mathbf{m}}

\newcommand\BasePt{w}
\newcommand\FiltPt{z}

\newcommand\StartPt{w}
\newcommand\EndPt{z}
\newcommand\StartPts{\mathbf{w}}
\newcommand\EndPts{\mathbf{z}}

\newcommand\is{\mathbf{i}}
\newcommand\js{\mathbf{j}}

\newcommand\bPsip{{}^b\Psi^+}
\newcommand\bPsia{{}^b{\widehat\Psi}}
\newcommand\bPsim{{}^b\Psi^-}
\newcommand\Psip{\Psi^+}
\newcommand\Psim{\Psi^-}

\hyphenation{ho-mol-o-gous}
\newcommand\commentable[1]{#1}
\newcommand\Rk{\mathrm{rk}}
\newcommand\Id{\mathrm{Id}}

\newcommand{\Tors}{\mathrm{Tors}}
\newcommand{\rk}{\mathrm{rk}}
\newcommand{\HF}{HF}
\newtheorem{theorem}{Theorem}
\newtheorem{thm}{Theorem}[section]
\newtheorem{prop}[thm]{Proposition}
\newtheorem{cor}[thm]{Corollary}
\newtheorem{conj}[thm]{Conjecture}
\newtheorem{lem}[thm]{Lemma}
\newtheorem{claim}[thm]{Claim}
\newtheorem{example}[thm]{Example}
\newtheorem{defn}[thm]{Definition}
\newtheorem{addendum}[thm]{Addendum}
\newtheorem{remark}[thm]{Remark}
\newtheorem{constr}[thm]{Construction}

\def\proof{\vspace{2ex}\noindent{\bf Proof.} }
\def\endproof{\relax\ifmmode\expandafter\endproofmath\else
  \unskip\nobreak\hfil\penalty50\hskip.75em\hbox{}\nobreak\hfil\bull
  {\parfillskip=0pt \finalhyphendemerits=0 \bigbreak}\fi}
\def\endproofmath$${\eqno\bull$$\bigbreak}
\def\bull{\vbox{\hrule\hbox{\vrule\kern3pt\vbox{\kern6pt}\kern3pt\vrule}\hrule}}
\newcommand{\smargin}[1]{\marginpar{\tiny{#1}}}

\newcounter{bean}
\newcommand{\Q}{\mathbb{Q}}
\newcommand{\R}{\mathbb{R}}
\newcommand{\T}{\mathbb{T}}
\newcommand{\Quat}{\mathbb{H}}
\newcommand{\C}{\mathbb{C}}
\newcommand{\N}{\mathbb{N}}
\newcommand{\Z}{\mathbb{Z}}

\newcommand{\OneHalf}{\frac{1}{2}}
\newcommand{\ThreeHalves}{\frac{3}{2}}
\newcommand{\OneQuarter}{\frac{1}{4}}
\newcommand{\CP}[1]{{\mathbb{CP}}^{#1}}
\newcommand{\CPbar}{{\overline{\mathbb{CP}}}^2}
\newcommand{\Zmod}[1]{\Z/{#1}\Z}

\newcommand{\Tr}{\mathrm{Tr}}
\newcommand{\Ker}{\mathrm{Ker}}
\newcommand{\CoKer}{\mathrm{Coker}}
\newcommand{\Coker}{\mathrm{Coker}}
\newcommand{\ind}{\mathrm{ind}}
\newcommand{\Image}{\mathrm{Im}}
\newcommand{\Span}{\mathrm{Span}}
\newcommand{\Spec}{\mathrm{Spec}}

\newcommand{\grad}{\vec\nabla}

\newcommand{\cm}{\cdot}
\newcommand\Sections{\mbox{$\Gamma$}}
\newcommand{\Nbd}[1]{{\mathrm{nd}}(#1)}
\newcommand{\nbd}[1]{\Nbd{#1}}
\newcommand{\CDisk}{D}

\newcommand{\Ideal}[1]{{\mathcal{I}}_{#1}}
\newcommand{\SheafRegFun}{{\mathcal O}}

\newcommand{\ModVort}{{\mathcal M}_{\mathrm{vort}}}
\newcommand{\ModSW}{\ModSWfour}
\newcommand{\ModSWirr}{\ModSW^{irr}}
\newcommand{\ModSWred}{\ModSW^{red}}
\newcommand{\ModSWthree}{\mathcal{N}}
\newcommand{\ModSWtwo}{N}
\newcommand{\ModSWfour}{\mathcal{M}}
\newcommand{\ModFlow}{\ModSWfour}
\newcommand{\SW}{SW}
\newcommand{\upm}{{\widehat {\ModFlow}}}
\newcommand{\HFswred}{HF^{SW} _{{\rm red}}}

\newcommand{\Sobol}[2]{L^{#1}_{#2}}

\newcommand{\SpinBunPinit}{\mbox{$\not\!S^+_0$}}
\newcommand{\SpinBunTwo}{\not\!S}
\newcommand{\SpinBunP}{\mbox{$\SpinBunFour^+$}}
\newcommand{\SpinBunM}{\mbox{$\SpinBunFour^-$}}
\newcommand{\SpinBunPM}{\SpinBunFour^\pm}
\newcommand{\SpinBun}{W}
\newcommand{\SpinBunFour}{W}
\newcommand{\SpinBunPFour}{\SpinBunFour^+}
\newcommand{\SpinBunPMFour}{\SpinBunFour^{\pm}}
\newcommand{\SpinBunMFour}{\SpinBunFour^-}

\newcommand{\Canon}[1]{{K}_{#1}}

\newcommand{\Dirac}{\mbox{$\not\!\!D$}}

\newcommand{\Conns}{\mathcal A}
\newcommand{\SWConfig}{\mbox{${\mathcal B}$}}
\newcommand{\SWConfigTwo}{\mbox{${\mathcal C}$}}
\newcommand{\Maps}{\mathrm{Map}}

\newcommand{\TrivBun}{\underline \C}
\newcommand{\TrivConn}{\underline d}
\newcommand{\SpinC}{{\mathrm{Spin}}^c}
\newcommand{\Spin}{{\mathrm{Spin}}}

\newcommand{\Proj}{\Pi}

\newcommand{\DDt}{\frac{\partial}{\partial t}}
\newcommand{\DDtheta}{\frac{\partial}{\partial\theta}}

\newcommand{\Deriv}{D}
\newcommand{\Vol}{\mathrm{Vol}}

\newcommand{\goesto}{\mapsto}

\newcommand{\Weitzenbock}{Weitzenb\"ock}
\newcommand{\Kahler}{K\"ahler}

\newcommand{\DBar}{\overline{\partial}}

\newcommand{\DDphi}{\frac{\partial}{\partial\phi}}

\newcommand{\CSD}{\mbox{${\rm CSD}$}}

\newcommand{\SpecFlow}{\mathrm{SF}}
\newcommand{\SF}{\SpecFlow}
\newcommand\Maslov{\mathrm{Mas}}
\newcommand\LinkingNumber{\mathrm{Lk}}
\newcommand\sgn{\mathrm{sgn}}

\newcommand\edim{\text{\rm{e-dim}}}

\newcommand\SWOp{\mathrm{sw}}
\newcommand\SWOpTwo{\mathrm{sw}}

\newcommand\Wedge{\Lambda}

\newcommand{\HessCSD}{D}
\newcommand\loc{\mathrm{loc}}
\newcommand\ext{\mathrm{ext}}
\newcommand\Pic[1]{\mathrm{Pic}^{#1}}
\newcommand\Hom{\mathrm{Hom}}

\newcommand{\SheafInvRegFun}[1]{\mbox{${\mathcal O}_{#1}^*$}}
\newcommand{\Cohom}[1]{\mbox{$H^{#1}$}}
\newcommand{\DerFun}[1]{\mbox{${\rm R}^{#1}$}}
\newcommand\abuts\Rightarrow
\newcommand\Sym{\mathrm{Sym}}

\newcommand\so{\mathfrak{so}}
\newcommand\Conn{\nabla}

\newcommand\End{\mathrm{End}}
\newcommand\Cliff{\rho}
\newcommand\Alg{\mathbb{A}}
\newcommand{\algel}{a}
\newcommand{\algelB}{b}
\newcommand{\Hol}{\mathrm{Hol}}

\newcommand\DDz{\frac{\partial}{\partial z}}
